\documentclass[12pt]{article}
\usepackage{latexsym,amssymb,amsmath,amsfonts,amsthm,amsopn,amscd}
\usepackage{graphicx} 
\newtheorem{lem}{Lemma}[section]
\newtheorem{thm}{Theorem}[section]

\newtheorem{defn}{Definition}[section]

\newcommand{\f}[1]{\mathfrak{#1}}

\newcommand{\up}{\widetilde{Sp}(n, \mb R)}
\newcommand{\mb}{\mathbb}
\newcommand{\commentout}[1]{}

\newcommand{\Hom}{\rm Hom \ }

\newcommand{\mc}{\mathcal}

\newcommand{\ot}{\otimes_{\w{Sp}(p, \mb R)}}

\newcommand{\Span}{{\rm span}}
\newcommand{\Ind}{{\rm Ind}}
\newcommand{\w}{\widetilde}
\begin{document}
\title{ Certain Induced Complementary Series of the Universal Covering of the Symplectic Group }
\author{Hongyu He 
\footnote{This research is supported in part by an NSF grant and LSU.}\\
Department of Mathematics,\\
 Louisiana State University, \\
 Baton Rouge, LA 70803, U.S.A.\\
email: livingstone@alum.mit.edu
}
\date{}
\maketitle 
\abstract{Let $\widetilde{Sp}(q, \mb R)$ be the universal covering of the symplectic group. In this paper, we study the unitarity problem for the representation induced from a one dimensional character $(\epsilon, v)$ of $\widetilde{GL}(q-p)$ tensoring with a unitary representation $\sigma_0$ of a smaller $\widetilde{Sp}(p, \mb R)$. We establish the unitarity when the real character $v$ is in a certain interval depending on $\epsilon$ and $\sigma_0$ satisfies a certain growth condition. We then apply our result inductively to construct complementary series for degenerate principal series with multiple $\w{GL}$-factors. In particular, in class $\epsilon$, there are $2^{q}$ principal complementary series of size {\it at least} $(0, c_{\epsilon})^{q}$ with $c_{\epsilon}=\min( |1- 2 \epsilon|, 1-|1-2\epsilon| )$. \\
\\
Various complementary series of the linear group $Sp(q, \mb R)$ have been constructed and studied by Kostant (~\cite{ko}), Knapp-Stein (~\cite{ks2}), Speh-Vogan (~\cite{sv}). Their complementary series are close to the tempered dual, while our complementary series are often \lq\lq far away\! \rq\rq from the tempered dual. More recently, Barbasch obtain all spherical unitary representations of $Sp(q, \mb R)$ (~\cite{barba}). Our approach is quite different and works well for the universal covering group. Essentially, we realize the underlying $\Ind$ functor as a Howe type duality with respect to a degenerate principal series representation $I(\epsilon, v)$ of $\w{Sp}(p+q, \mb R)$ (~\cite{howe}). Then we construct an induced intertwining operator for the induced representation under consideration from the intertwining operator on $I(\epsilon, v)$. The positivity of the induced intertwining operator is established by a standard deformation argument based on the positivity of the intertwining operator on $I(\epsilon, v)$. }

\section{Introduction} 
For any semisimple Lie group $G$, let $\Pi(G)$ be the admissible dual, $\Pi_u(G)$ be the unitary dual and $\Pi_2(G)$ be the tempered dual. Let $Sp$ be the real symplectic group. Let $\widetilde{Sp}$ be the universal covering of $Sp$. For any subgroup $H$ of $Sp$, let $\tilde H$ be the preimage under the universal covering. Let $\mc C$ be the preimage of the identity. Then $\mc C \cong \mathbb Z$. 
We parametrize the unitary dual of $\mathbb Z$ by $\epsilon \in [0, 1)$. Then $\Pi_{*}(\w{Sp})$ is a union of
$\Pi_{*}(\w{Sp})_{\epsilon}$ on which $\mc C$ acts by $\epsilon$. For instance, $\Pi_{*}(\w{Sp})_{0} = \Pi_{*}(Sp)$ and $\Pi_{*}(\w{Sp})_{\frac{1}{2}}$ is the genuine dual of the metaplectic groups. If $\mc C$ acts on a representation by $\epsilon$, we say that this representation is of {\bf class $\epsilon$}. \\
\\
Let $P$ be a parabolic subgroup of $Sp(n, \mb R)$. Let $MAN$ be the Langlands decomposition of $P$. Let $\sigma \in \Pi_u(\tilde M)$. Let $v \in \f a^*_{\mathbb C}$. Let $\Ind_{\tilde P}^{\w{Sp}(n, \mb R)} \sigma \otimes \mathbb C_{v}$ be the normalized induced representation. When $v \in i \f a^*$, $\Ind_{\tilde P}^{\w{Sp}(n, \mb R)} \sigma \otimes \mathbb C_{v}$ is always unitary. This induction process is called unitary parabolic induction. The main theme of this paper is to discuss the unitarity of $\Ind_{\tilde P}^{\w{Sp}(n, \mb R)} \sigma \otimes \mathbb C_{v}$ for $v$ real. Such a series of unitary representations is often called complementary series. I should point out that the nonunitary $\Ind_{\tilde P}^{\w{Sp}(n, \mb R)} \sigma \otimes \mathbb C_{v}$ may contain unitarizable subquotients. We are not concerned with the unitarizable subquotients in this paper.\\
\\
Complementary series for $SL(2, \mb R)$ were first constructed by V. Bargmann and Gelfand-Naimark. For the universal covering of $SL(2, \mb R)$, a series of length $|1-2 \epsilon |$ was found in each $\Pi_u(\w{SL}(2, \mb R))_{\epsilon}$ by Puk\'anszky (~\cite{puk}). These complementary series are nondegenerate deformations of the unitary principal series. However, there are also complementary series that are not nondegenerate deformation of the unitary principal series, as pointed out by G. Zuckerman. In ~\cite{du}, M. Duflo found two unitarizable regions in the principal series of the complex $G_2$. One cannot deform the Hilbert structure to go from one region to the other without degeneracy. The same phenomena occurs for real $G_2$ (~\cite{vogang2}). It is not clear whether this phenomena occurs for classical groups. At least, for general linear groups, all complementary series are deformations of unitarily induced representations, according to Vogan's classification of unitary dual (~\cite{vogan}). \\
\\
Systematical studies of complementary series for semisimple Lie groups came about in a series paper by Kostant (~\cite{ko}), Knapp-Stein (~\cite{ks}, ~\cite{ks2}) and Speh-Vogan (~\cite{sv}).
They all obtained powerful results about the complementary series, mainly near the tempered dual. The other extreme is the degenerate complementary series, which was under intensive investigation by many authors. See  ~\cite{kr}, ~\cite{john}, ~\cite{bss}, ~\cite{oz}, ~\cite{boo}, ~\cite{lee}, ~\cite{lz}, ~\cite{hl} and references therein.  Perhaps, one of the most well-known examples is Stein's complementary series, which are induced from a one dimensional real character of two copies of $GL(n)$. For $\w{Sp}$, Sahi gives a complete classification of complementary series induced from a one dimensional character of the Siegel parabolic subgroup, which in a way, resembles 
Stein representations (~\cite{sahi}). \\
\\
In this paper, we study the following induced representation
$$\Ind_{\w{Sp}(p, \mb R)\w{GL}(q-p)N}^{\w{Sp}(q, \mb R)} \sigma \otimes \mu^{\epsilon} \otimes {\nu}^{t},$$
where $\sigma \in \Pi_u(\w{Sp}(p, \mb R))$, $\mu^{\epsilon}$ is a character on the component group of $\w{GL}(q-p)$ that is compatible with $\sigma$, and $\nu$ is the $|\det|$ on $\w{GL}(q-p)$. {\bf $\mu^{\epsilon}$ being compatible with $\sigma$ means that
$\mu^{\epsilon}|_{\mc C}=\sigma|_{\mc C}$ as scalars}. Before we state our main result, we fix some notations. Let $\mathbf a$ be the constant vector $(a, a, \ldots, a)$ of a suitable size. For two $n$-dimensional real vectors $\lambda$ and $\mu$, we say that
$\lambda \prec \mu$ if $\sum_{i=1}^{l} \lambda_i < \sum_{i=1}^{l} \mu_i$ for every $l \in [1, n]$. We say that $\lambda \preceq \mu$ if $\sum_{i=1}^{l} \lambda_i \leq \sum_{i=1}^{l} \mu_i$ for every $l \in [1, n]$. 

\begin{thm} Let $\epsilon \in [0, 1)$ and 
Let $\sigma \in \Pi_u(\w{Sp}(p, \mb R))$. Suppose that $\sigma$ is compatible with $\mu^{\epsilon}$ and every leading exponent $v$ of $\sigma$ satisfies
$$\Re (v)-\mathbf{\frac{p+q+1}{2}}+ (2p, 2p-2, \ldots 2) \prec 0.$$
Suppose that $|t| <  \frac{1}{2}-|\frac{1}{2}-|2 \epsilon-1|| $ if $q-p$ is odd or $|t| < |\frac{1}{2}-|2 \epsilon-1||$ if $q-p$ is even. Then $\Ind_{\w{Sp}(p, \mb R)\w{GL}(q-p)N}^{\w{Sp}(q, \mb R)} \sigma \otimes \mu^{\epsilon} \otimes {\nu}^{t}$ is unitarizable. 
\end{thm}
I shall point out that our condition on $t$ is only a sufficient condition. The complementary series can be longer than the one given in our theorem. \\
\\
If $\sigma$ is tempered, its leading exponents always satisfy the inequalities specified in our theorem. 
Now one can apply our theorem inductively to obtain complementary series induced from 1 dimensional characters of the $\w{GL}$-factors.
I shall mention one special case with only $GL(1)$-factors. The general statement is given in Theorem ~\ref{mu1}.
\begin{thm}  Let $\sigma$ be a tempered representation of $\w{Sp}(n-r, \mb R)$. Suppose that the characters ${\mu}^{\epsilon_i}$ of $\w{GL}(1)$ are compatible with $\sigma$. Then
$$\Ind_{\w{Sp}(n-r, \mb R) \w{GL}(1)^r N}^{\w{Sp}(n, \mb R)} \sigma \otimes ( \otimes_{i=1}^r {\mu_i}^{\epsilon_i} \otimes {\nu_i}^{t_i})$$
is unitarizable for all $t_i \in (-c_{\epsilon_i}, c_{\epsilon_i})$ with $c_{\epsilon_i}= \min( |1- 2 \epsilon_i|, 1-|1-2\epsilon_i| )$. 
\end{thm}
The representations in this example are all close to the tempered dual and the results for finite coverings are known (~\cite{ks2}, ~\cite{sv}). If we only look at the metaplectic group, there are $2^{r}$ genuine complementary series of size $(0, \frac{1}{2})^r$. These complementary series seem to be all there are for this particular induction. For the linear group $Sp(n, \mb R)$, our result says nothing about the complementary series. It remains an interesting problem to investigate why this happens.\\
\\
The following is what is covered in this paper. Let $p< q$ and $p+q=n$. In Section 2, we discuss the degenerate principal series $I(\epsilon, t)$ of $\w{Sp}(n, \mb R)$ and its restriction to $\w{Sp}(p, \mb R) \w{Sp}(q, \mb R)$ (~\cite{comp}). In particular, we review some results from ~\cite{comp} regarding the mixed model on which the action of $\w{Sp}(p, \mb R) \w{Sp}(q, \mb R)$ is explicit (See Theorem \!~\ref{mix}). In Section 3, we discuss the parabolic induction for $\w{Sp}(q,\mb R)$, whose center is infinite. We then decompose the unitary degenerate principal series $I(\epsilon, i t)$ into a direct integral based on the action of $\w{Sp}(p, \mb R) \w{Sp}(q, \mb R)$ (See Theorem \!~\ref{pla}). This gives an $L^2$ Howe Type duality, which can be identified with a unitary parabolic induction. In Section 4, we construct the invariant tensor functor
$$ \sigma \in \Pi_u(\w{Sp}(p, \mb R)) \rightarrow I^{\infty}(\epsilon, t) \ot V(\sigma)$$
under {\bf a growth condition} on $\sigma$ (see Definition ~\ref{itp}). We prove that 
$$I^{\infty}(\epsilon, t)_{\w{U}(q)} \ot V(\sigma) \cong [\Ind_{\w{Sp}(p, \mb R)\w{GL}(q-p)N}^{\w{Sp}(q, \mb R)} \pi \otimes \mu^{\epsilon} \otimes {\nu}^{t})]_{\w{U}(q)}$$
(see Theorem \!~\ref{id}). In Section 5, we define an induced intertwining operator 
$$A(\epsilon, t, \sigma): I^{\infty}(\epsilon, t)_{\w{U}(q)} \ot V(\sigma) \rightarrow I^{\infty}(\epsilon, -t)_{\w{U}(q)} \ot V(\sigma)$$
from $A(\epsilon, t): I^{\infty}(\epsilon, t)_{\w{U}(q)} \rightarrow I^{\infty}(\epsilon, -t)_{\w{U}(q)}$. We show that $A(\epsilon, t, \sigma)$ inherits the positivity of $A(\epsilon, t)$. Hence, $\Ind_{\w{Sp}(p, \mb R)\w{GL}(q-p)N}^{\w{Sp}(q, \mb R)} \pi \otimes \mu^{\epsilon} \otimes {\nu}^{t}$ is unitarizable if $I(\epsilon, t)$ is unitarizable (See Theorem \!~\ref{main0}). In Section 6, we define a more general version of
$A(\epsilon, t, \sigma)$ for all $t$ sufficiently negative and sharpen our results from Section 5 (See Theorem \!~\ref{main}). How $A(\epsilon, t, \sigma)$ fits into the general theory of intertwining operator remains to be an interesting problem. In Section 7, we construct some induced complementary series inductively. See Theorem \!~\ref{mu1}. \\
\\
This paper is based on a talk given at Yale. I would like to thank Prof. Roger Howe and Prof. Gregg Zuckerman for their interests.

\section{Degenerate Complementary Series}
Let $\mathbb C^n$ be the $n$-dimensional complex Hilbert space with an orthonormal basis
$$\{e_1, e_2, \ldots, e_n \}.$$
Let $\Omega(\, , \,)=\Im(\, , \,).$
Regarding $\mathbb C^n$ as a real vector space, let $Sp(n, \mb R)$ be the symplectic group that preserves $\Omega$. Let $U(n)$ be the unitary group that preserves $(\, , \,)$. Clearly $U(n)$ is a maximal compact subgroup of $Sp(n, \mb R)$. \\
\\
Let $P$ be the Siegel parabolic subgroup that preserves the real linear span of
$$\{  i e_1, i e_2, \ldots, i e_n \},$$
and let $N$ be its nilradical. Choose the Levi factor to be the subgroup of $P$ that preserves the real linear span of $\{ e_i \mid i \in [1, n] \}$. Clearly, $L \cong GL(n, \mb R)$ and $L \cap U(n) \cong O(n)$. \\
 \\
 On the covering group, we have $\tilde L \cap \tilde{U}(n)=\tilde{O}(n)$. Recall that
$$\tilde{U}(n)= \{ (x, g) \mid g \in U(n), \exp 2 \pi i x = \det g, x \in \mathbb R \}.$$
Therefore
$$\tilde{O}(n)=\{ (x, g) \mid g \in O(n), \exp 2 \pi i x= \det g, x \in \mathbb R \}.$$
Notice that for $g \in O(n)$, $\det g= \pm 1$. So $ x \in \frac{1}{2} \mathbb Z$. Identify the identity component of $\tilde O(n)$ with $SO(n)$.
We have the following exact sequence
$$ 1 \rightarrow SO(n) \rightarrow \tilde{O}(n) \rightarrow \frac{1}{2} \mathbb Z \rightarrow 1.$$
Consequently, we have
$$ 1 \rightarrow GL_0(n, \mb R) \rightarrow \tilde L \rightarrow \frac{1}{2} \mathbb Z \rightarrow 1.$$
In fact, $$\tilde L=\{ (x, g) \mid g \in L, \exp 2 \pi i x= \frac{\det g}{|\det g|}, x \in \mathbb R \}.$$
The one dimensional unitary characters of $\frac{1}{2} \mathbb Z$ are parametrized by the one dimensional torus $T$. Identify $T$ with $[0, 1)$. {\bf Let $\mu^{\epsilon}$ be the character of $\frac{1}{2} \mathbb Z$ corresponding to $\epsilon \in [0,1)$}.
Now each character $\mu^{\epsilon}$ yields a character of $\tilde L$, which in turn, yields a character of $\tilde P$. For simplicity, we retain $\mu^{\epsilon}$ to denote the character on $\tilde L$ and $\tilde P$. Let $\nu$ be the $\det$-character on ${\tilde L}$, i.e.,
\begin{equation}~\label{deter}
\nu(x,g)= | \det g | \qquad (x,g) \in \tilde L.
\end{equation}
Let
$I(\epsilon, t)=\Ind_{\tilde P}^{\up} \mu^{\epsilon} \otimes \nu^{t}$
be the normalized induced representation of $\w{Sp}(n, \mb R)$ with $\epsilon \in [0,1)$ and $t \in \mathbb C$ (~\cite{sahi}).  $I(\epsilon,t)$ is called a degenerate principal series representation. 
Clearly, $I(\epsilon, t)$ is unitary when $t \in i \mathbb R$.

\begin{thm}[Thm A, ~\cite{sahi}]
Suppose that $t$ is real.
For $n$ even, $I(\epsilon,t)$ is irreducible and unitarizable if and only if 
$0< |t| <|\frac{1}{2}-|2 \epsilon-1||$. For $n$ odd and $n >1$, $I(\epsilon, t)$ is irreducible and unitarizable if and only if $0< |t| < \frac{1}{2}-|\frac{1}{2}-|2 \epsilon-1||$.
\end{thm}
See also 
 ~\cite{kr}, ~\cite{boo}, ~\cite{oz}, ~\cite{lee} and the references therein.

\begin{defn}~\label{dia} Let $n=p+q$. Let $\mathbb C^p$ be the complex linear space spanned by $\{e_1, e_2, \ldots, e_p \}$ and $\mathbb C^q$ be the complex linear space spanned by $\{e_{p+1}, e_{p+2}, \ldots, e_{n} \}$. Let $\Omega_p=-\Omega|_{\mathbb C^p}$ and $\Omega_q=\Omega|_{\mathbb C^q}$. Write
$$\Omega=-\Omega_p+\Omega_q.$$
 Let $Sp(p, \mb R)$ be the symplectic group preserving $-\Omega_p$ and fixing $\mathbb C^q$;
$Sp(q, \mb R)$ be the symplectic group preserving $\Omega_q$ and fixing $\mathbb C^p$.
 We say that $(Sp(p, \mb R), Sp(q, \mb R))$ is {\it diagonally embedded } in $Sp(n, \mb R)$. 
 Let $U(p)=U(n) \cap Sp(p, \mb R)$ and $U(q)=U(n) \cap Sp(q, \mb R)$.
\end{defn}
Although the symplectic group preserving $\Omega_p$ also preserves $-\Omega_p$, the parametrization of $Sp(-\Omega_p)$ will be according to the bases $\{ e_1,e_2, \ldots e_p, i e_1, i e_2, \ldots i e_p \}$, not the standard basis $\{ i e_1, i e_2, \ldots, i e_p, e_1, e_2, \ldots, e_p \}$. The reader should note that this difference of parametrization will incur an involution on the representation level if we stick with the standard basis (see ~\cite{nu}).\\
\\
{\bf Suppose from now on that $p < q$ and $p+q=n$}. Some of the statements do make sense for $p=q$.
Consider the action of $\w{Sp}(p, \mb R) \w{Sp}(q, \mathbb R)$ on $I(\epsilon, t)$. Recall that
$I^{\infty}(\epsilon, t)$ consists of smooth sections of
the homogeneous line bundle $\mc L_{\epsilon, t}$
$$ \up \times_{\tilde P} \mathbb C_{\mu^{\epsilon} \otimes \nu^{t+\rho}} \rightarrow X,$$
where $\rho=\frac{n+1}{2}$. Identify $X$ with the variety of Lagrangian Grassmanian. Choose a base point 
$$x_0=\Span_{\mathbb R} \{i \, e_j+e_{p+j},  e_j + i \, e_{p+j}, i e_k \mid j \in [1, p]; k \in [2p+1, n] \}.$$  
Let $P_{q-p}(q)$ be the maximal parabolic subgroup of $Sp(q, \mathbb R)$ that preserves
$$W= \Span_{\mathbb R} \{ i e_k \mid k \in [2p+1, n] \}.$$
Let $W^{\perp}=\Span_{\mathbb R} \{ e_j, i e_j, i e_k \mid j \in [p+1, 2p], k \in [2p+1, k] \}$.
Let $Q_{q-p}(q)$ be the subgroup of $P_{q-p}(q)$ that fixes every vector in $W^{\perp}/W$. Let $N_{q-p}(q)$ be the nilradical of $P_{q-p}(q)$. Then the stabilizer
$$Sp(q, \mathbb R)_{x_0}=Q_{q-p}(q).$$
In addition, identifying the basis $i e_j$ with $e_{j+p}$ and $e_j$ with $ i e_{j+p}$ for every $j \in [1,p]$, we obtain a
symplectic isomorphism from $(\mathbb C^{p}, -\Omega_p)$ onto a linear subspace of $(\mathbb C^{q}, \Omega_q)$. This isomorphism induces a group isomorphism:
$$ g \in Sp(p, \mathbb R)   \rightarrow \dot{g} \in Sp(q, \mb R).$$
Let $\Delta(Sp(p, \mathbb R))=\{(g, \dot g) \mid g \in Sp(p, \mathbb R) \} \subset Sp(p, \mb R) Sp(q, \mb R).$ 
Then the isotropy group 
$$(Sp(p, \mb R) Sp(q, \mb R))_{x_0}= \Delta(Sp(p, \mb R)) Q_{q-p}(q).$$
See ~\cite{nu} ~\cite{comp}. \\
\\
Take the Siegel parabolic subgroup $P$ to be the stabilizer of $x_0$. Then
$$P \cap Sp(p, \mb R) Sp(q, \mb R) = \Delta(Sp(p, \mb R)) Q_{q-p}(q).$$
We make two observations regarding the covering. First of all, since $\Delta(Sp(p, \mb R))$ preserves any top degree exterior products in the Lagrangian $x_0$, we have
$$\nu(\Delta(Sp(p, \mathbb R)))= 1.$$
 Hence
$\w{\Delta(Sp(p, \mb R))}$ splits into a direct product
$\mc C \times \w{\Delta(Sp(p, \mb R))}_0$
where $\w{\Delta(Sp(p, \mb R))}_0$ is the identity component of $\w{\Delta(Sp(p, \mb R))}$, which can be identified with $\Delta(Sp(p, \mb R))$. Secondly, the group $Q_{q-p}(q)$ has a Levi decomposition $GL(q-p) N_{q-p}(q)$. It has two connected components contained in the two connected components of $P$ respectively. It follows that $\nu(h)$ coincides with $|\det(h)|$ for any $h \in \w{GL}(q-p) \subset \w{Q}_{q-p}(p)$. \\
\\
Now we can restrict the line bundle $\mc L_{\epsilon, t}$ onto $\w{Sp}(q, \mb R)$. We obtain
\begin{equation}~\label{mt}
\mc M_{\epsilon, t}: \widetilde{Sp}(q, \mb R) \times_{\widetilde{GL}(q-p) N_{q-p}(q)} \mathbb C_{\mu^{\epsilon} \otimes \nu^{t+\rho}}
\rightarrow Sp(q, \mb R)/Q_{q-p}(q).
\end{equation}
Notice that $Sp(q, \mb R)/Q_{q-p}(q)$ has a principal bundle structure
$$Sp(p, \mb R) \rightarrow Sp(q, \mb R)/Q_{q-p}(q) \rightarrow Sp(q, \mb R)/P_{q-p}(q) \cong U(q)/U(p)O(q-p).$$
We parametrize $Sp(q, \mb R)/Q_{q-p}(q)$ by 
$$ [ g_1 \in Sp(p, \mb R), k_2 \in U(q) ] \in (Sp(p, \mb R), U(q)/O(q-p))/U(p) ,$$
and equip it with the invariant measure $ d g_1 d [k_2] $.  The parametrization is given by
$$ [g_1, k_2] \rightarrow k_2 g_1 Q_{q-p}(q)$$
since $Sp(p, \mb R) \cap U(q) = U(p)$.
\begin{thm}[Page 11-12 ~\cite{comp}]~\label{mix} 
Let $t \in \mb R$, $p<q$ and $p+q=n$. 
\begin{enumerate}
\item The restriction map $f \rightarrow f|_{\w{Sp}(q, \mb R)}$ induces an isometry between $I(\epsilon, i t)$  and $ L^2(\mc M_{\epsilon, it}, d [g_1] d [k_2])$. Let $\w{Sp}(q, \mb R)$ act on  $ L^2(\mc M_{\epsilon, it}, d [g_1] d [k_2])$ from the left and let $\w{Sp}(p, \mb R)$ act on $ L^2(\mc M_{\epsilon, i t}, d [g_1] d [k_2])$ from the right. Then the restriction map intertwines the actions of $\widetilde{Sp}(p, \mb R) \widetilde{Sp}(q, \mb R)$. So as $\widetilde{Sp}(p, \mb R) \widetilde{Sp}(q, \mb R)$ representations, $I(\epsilon, i t) \cong  L^2(\mc M_{\epsilon, it}, d [g_1] d [k_2]).$ 
\item Let $\w{Sp}(q, \mb R)$ act on the space of smooth sections $ C^{\infty}(\mc M_{\epsilon, t}, d [g_1] d [k_2])$ from the left and let $\w{Sp}(p, \mb R)$ act on $ C^{\infty}(\mc M_{\epsilon,  t}, d [g_1] d [k_2])$ from the right.
Then for every $f \in I^{\infty}(\epsilon, t)$, $f|_{\w{Sp}(q, \mb R)}$ is smooth and $|f(g_1 k_2)|$ is bounded by a multiple of
$$\det(g_1 g_1^t+ I)^{-\frac{n+1+t}{4}}.$$
The restriction map $f \rightarrow f|_{\w{Sp}(q, \mb R)}$ intertwines $I^{\infty}(\epsilon, t)$ with 
$C^{\infty}(\mc M_{\epsilon, t}, d [g_1] d [k_2])$.
\end{enumerate}
\end{thm}
 The restriction map in (2) is not onto. Its image is the space of smooth functions satisfies the decaying condition
specified in (2) and a certain analytic condition at $\infty$. We call the model on $\mc M_{\epsilon, t}$ {\bf mixed model}.

\section{Parabolic Induction for $\widetilde{Sp}$ and a Howe Type Duality}
Let $P$ be a parabolic subgroup of $Sp(n, \mb R)$.  Then the Levi subgroup $L$ is of the form
$$GL(r_1) \times GL(r_2) \times \ldots \times GL(r_s) \times Sp(r_0, \mb R)$$
with $\sum_{i=1}^{s} r_i =n-r_0$ and $r_i \geq 1$. 
The covering $\tilde L$ is of the form
$$\w{GL}(r_1) \times_{\mc C} \w{GL}(r_2) \times_{\mc C} \ldots \times_{\mc C} \w{GL}(r_s) \times_{\mc C} \w{Sp}(r_0, \mb R);$$
$$\tilde P \cong \tilde L N.$$
Here
$H_1 \times_{\mc C} H_2$ is defined to be the quotient group
$$ H_1 \times H_2 / \{ (c, c^{-1}) \mid c \in \mc C \}$$
whenever the subgroup $\mc C$ is in the center of both $H_1$ and $H_2$. \\
\\
Now an irreducible admissible representation of $\tilde L$ is of the following form:
$$\sigma_1 \otimes_{\mc C} \sigma_2 \ldots \otimes_{\mc C} \sigma_s \otimes_{\mc C} \sigma_0,$$
where $\sigma_i|_{\mc C}$ and $\sigma_0|_{\mc C}$ are all scalar multiplications and the underlying one dimensional characters are all the same. In this situation, we say that $\{\sigma_i, \sigma_0 \}$ are compatible. For simplicity, we will write
$$\sigma_1 \otimes \sigma_2 \ldots \otimes \sigma_s \otimes \sigma_0,$$
for $\sigma_1 \otimes_{\mc C} \sigma_2 \ldots \otimes_{\mc C} \sigma_s \otimes_{\mc C} \sigma_0,$
and we will {\bf always assume that $\{\sigma_i, \sigma_0 \}$ are compatible}.\\
\\
Now let $\{\sigma_i, \sigma_0\}$ be compatible Hilbert representations. We can define parabolic induction
$$I(\sigma_i, \sigma_0)=\Ind_{\tilde P}^{\w{Sp}(n, \mb R)} [\otimes_{i} \sigma_i] \otimes \sigma_0.$$
The smooth vectors $I^{\infty}(\sigma_i, \sigma_0)$ consists of smooth sections of
$$\w{Sp}(n, \mb R) \times_{\tilde P} [\otimes_{i} \sigma_i] \otimes \sigma_0 \rightarrow \w{Sp}(n, \mb R)/\tilde P.$$
Let us analyze the one dimensional characters of $\w{GL}(r_i)$. Notice that $\mc C \subseteq \w{GL}(r_i)$. If one identifies $\mc C$ with $\mathbb Z$, then the component group of $\w{GL}(r_i)$ can be identified with $\frac{1}{2} \mathbb Z$. Now parametrize the unitary dual of
$\w{GL}(r_i)/\w{GL}(r_i)_0$ by a real number $\epsilon \in [0, 1)$. As in the Siegel parabolic case, let $\mu^{\epsilon}$ be the character of $\frac{1}{2} \mathbb Z$ corresponding to $\epsilon \in [0,1)$. 
Two one dimensional character $\mu^{\epsilon_1}$ and $\mu^{\epsilon_2}$ are compatible if and only if $ 2 \epsilon_1- 2 \epsilon_2$ is an integer. For any $\sigma_0$ of $\w{Sp}(r_0, \mb R)$, there are only {\bf two} $\epsilon_i$ such that ${\mu}^{\epsilon_i}$ is compatible with $\sigma_0$.\\
\\
We shall now define a {\bf universal} $\mu$ for $\w{GL}$ that are compatible with any $\w{GL}(r) \hookrightarrow \w{GL}(s) ( r \leq t) $. Now each character $\prod_{i} {\mu_i}^{\epsilon_i}$ yields a character of $\tilde L$, which in turn, yields a character of $\tilde P$. For simplicity, we retain $\prod_{i}{\mu_i}^{\epsilon_i}$ to denote the character on $\tilde L$ and $\tilde P$. Let $\nu$ be the {\bf universal} $|\det|$-character for $\w{GL}$. For $t_i \in \mathbb C$, we can define a character $\prod_i \nu_i^{t_i}$ for $\tilde L$ and $\tilde P$. \\
\\
We may now decompose $I(\epsilon, i t)$ as a direct integral of irreducible unitary representations of $\w{Sp}(p, \mb R) \w{Sp}(q, \mb R)$:
\begin{thm}~\label{pla} 
Suppose that $p+q=n$ and $q>p$. Then
$$I(\epsilon, i t)|_{\w{Sp}(p, \mb R) \w{Sp}(q, \mb R)} \cong \int_{\sigma \, tempered}  [\Ind_{\w{Sp}(p, \mb R)\w{GL}(q-p)N}^{\w{Sp}(q, \mb R)} \sigma^* \otimes {\mu}^{\epsilon} \otimes {\nu}^{it} ] \otimes \sigma d \sigma.$$
Here $\sigma^*|_{\mc C}={\mu}^{\epsilon}|_{\mc C}$ as scalars and
$d \sigma$ is a certain measure on $\Pi_2(\w{Sp}(p, \mb R))$.
\end{thm}
This Theorem can be easily derived from Theorem ~\ref{mix} and Equ. ~\ref{mt}. \\
\\
Now the restriction $I(\epsilon, i t)_{\w{Sp}(p, \mb R) \w{Sp}(q, \mb R)}$ sets up an $L^2$-version Howe type duality:
$$\sigma \in \Pi_2(\w{Sp}(p, \mb R)) \rightarrow \Ind_{\w{Sp}(p, \mb R)\w{GL}(q-p)N}^{\w{Sp}(q, \mb R)} \sigma^* \otimes {\mu}^{\epsilon} \otimes {\nu}^{it} \in \Pi_u(\w{Sp}(q, \mb R)).$$
Two problems arise. The first is to extend this correspondence to representations beyond tempered ones. One could, for instance, define a correspondence algebraically, like in ~\cite{howe}. In our situation, this will be tautological.  The second is a Plancherel formula for the degenerate complementary series $I(\epsilon, t)$. One would expect a correspondence of the following form:
$$\sigma \in \Pi_2(\w{Sp}(p, \mb R)) \rightarrow \Ind_{\w{Sp}(p, \mb R)\w{GL}(q-p)N}^{\w{Sp}(q, \mb R)} \sigma^* \otimes {\mu}^{\epsilon} \otimes {\nu}^{t} \in \Pi_u(\w{Sp}(q, \mb R)).$$
However, to prove a statement like this would take us on a completely different route. To stick with our main goal, that is, to prove the unitarity of the right hand side for certain $\sigma$ and $t$, we will study the following invariant tensor functor
$$\sigma \rightarrow I^{\infty}(\epsilon, t) \otimes_{\w{Sp}(p, \mb R)} V(\sigma)$$
which we will define in the next section. This scheme of construction is motivated by ~\cite{li2}, ~\cite{theta}, ~\cite{unit}.

\section{Invariant Tensor Product and Parabolic Induction}
Given an irreducible admissible representation $\pi$ of $\w{Sp}(p, \mb R)$, let $V(\pi)$ be the space of $\w{U}(p)$-finite vectors. $V(\pi)$ is a direct sum of irreducible $\w{U}(p)$-modules with finite multiplicities. It is often called the Harish-Chandra module of $\pi$. Let $I^{\infty}(\epsilon, t)$ be the smooth vectors in the mixed model of $I(\epsilon, t)$. 
\begin{defn}~\label{itp}
Let $f \in I^{\infty}(\epsilon, t)$. Let $\pi \in \Pi_u(\w{Sp}(p, \mb R))$ be compatible with $\mu^{\epsilon}$. Let $v, u \in V(\pi)$. Suppose that 
$$\int_{\w{Sp}(p, \mb R)/\mc C} | \det(g_1 g_1^t+ I)^{-\frac{n+1+t}{4}} (\pi(g)v, u)| d [g] < \infty.$$
Then we define
$f \otimes_{\w{Sp}(p, \mb R)} v$ to be a $\Hom_{\mb C}(\pi^*, \mathbb C)$-valued function on $\w{Sp}(q, \mb R)$ as follows:
$$((f \otimes_{\w{Sp}(p, \mb R)} v)(g), u)=\int_{\w{Sp}(p, \mb R)/\mc C} f(g h)(\pi(h)v, u) d [h].$$
Here $\pi^*$ is the contragredient representation in the category of Harish-Chandra modules and $(\, , \,)$ is the pairing between $\pi$ and $\pi^*$. Since $\pi$ is unitary, $(\, , \,)$ can be identified with the inner product on $V(\pi)$. Let
$I^{\infty}(\epsilon, t) \otimes_{\w{Sp}(p, \mb R)} V(\pi)$ be the space spanned by $f \otimes_{\w{Sp}(p, \mb R)} v$.
\end{defn}
Let me make one  remark here. If $h \in \mc C$, we see that
$$f(g h)(\pi(h)v, u)=\mu^{\epsilon}(h^{-1}) f(g) \mu^{\epsilon}(h)(v, u)=f(g)(v, u).$$
Hence the integral $\int_{\w{Sp}(p, \mb R)/\mc C} f(g h)(\pi(h)v, u) d [h]$ is well-defined. \\
\\
Let me give a necessary condition and a sufficient condition for the well-definedness of $I^{\infty}(\epsilon, t) \ot V(\pi)$.
See Ch. 8.12 and 8.13 ~\cite{kn}.  
\begin{lem}~\label{itp1} Let $\pi \in \Pi_u(\w{Sp}(p, \mb R))$. Suppose that every leading exponent $v$ of $\pi$ satisfies
$$\Re (v)-\mathbf{\frac{n+1+t}{2}}+ 2 \rho(p) \prec 0.$$
Then $I^{\infty}(\epsilon, t) \ot V(\pi)$ is well-defined. Conversely, if $I^{\infty}(\epsilon, t) \ot V(\pi)$ is well-defined,
 then we have $$\Re (v) -\mathbf{\frac{n+1+t}{2}}+ 2\rho(p)  \preceq 0 .$$
 Here $\rho(p)=(p, p-1, \ldots, 1)$.
\end{lem}
If $\pi$ is tempered, then the leading exponents of $\pi$ are all bounded by $-\rho(p)$. In addition, if  $I(\epsilon, t)$ is unitarizable, then 
$$-\mathbf{\frac{n+1 \pm t}{2}} \preceq -\mathbf {\frac{2p+2 \pm t}{2}} \prec -\rho(p).$$
So for tempered $\pi$ and unitarizable $I(\epsilon, t)$, $I^{\infty}(\epsilon, t) \ot V(\pi)$ is always well-defined. \\
\\
{\bf From now on, we will assume that $I^{\infty}(\epsilon, t) \otimes_{\w{Sp}(p, \mb R)} V(\pi)$ is well-defined.} 
$f \otimes_{\w{Sp}(p, \mb R)} v$ is a function
$$ g \in \w{Sp}(q, \mb R) \rightarrow \int_{\w{Sp}(p, \mb R)/\mc C} f(g h) \pi(h)v d [h]$$
with value in $\Hom_{\mb C}(\pi^*, \mathbb C)$. In view of Equ. \!~\ref{mt}, we have the following lemma.
\begin{lem} Let $f \in I^{\infty}(\epsilon, t)$ and $v \in V(\pi)$. Then $$f \otimes_{\w{Sp}(p, \mb R)} v \in \Ind_{\w{Sp}(p, \mb R)\w{GL}(q-p)N}^{\w{Sp}(q, \mb R)} 
\Hom_{\mb C}(\pi^*, \mathbb C) \otimes \mu^{\epsilon} \otimes {\nu}^{t}.$$
\end{lem}
In this lemma, $\Ind$ means the space of sections because $\Hom_{\mb C}(\pi^*, \mathbb C)$ does not have a Hilbert space structure.
\begin{lem}
Suppose that $f \in I^{\infty}(\epsilon, t)_{\w{U}(q)}$, that is, $f$ is smooth and $\w{U}(q)$-finite. Suppose that $I^{\infty}(\epsilon, t) \otimes_{\w{Sp}(p, \mb R)} V(\pi)$ is well-defined. Then 
$f \otimes_{\w{Sp}(\mb R)} v$ is $\tilde U(q)$-finite. So $f \otimes_{\w{Sp}(\mb R)} v$ is in the Harish-Chandra module
$$ V(\Ind_{\w{Sp}(p, \mb R)\w{GL}(q-p)N}^{\w{Sp}(q, \mb R)} \pi \otimes \mu^{\epsilon} \otimes {\nu}^{t}).$$
\end{lem}
Proof: It is easy to see that the map $f \otimes v \rightarrow f \otimes_{\w{Sp}(p, \mb R)} v$ preserves the action of $\w{U}(q)$. If $f$ is $\w{U}(q)$-finite, then $f \otimes v$ is  $\w{U}(q)$-finite. Therefore $f \otimes_{\w{Sp}(p, \mb R)} v$ is also $\tilde U(q)$-finite.
The $\tilde U(q)$-finite vectors in 
$$\Ind_{\w{Sp}(p, \mb R)\w{GL}(q-p)N}^{\w{Sp}(q, \mb R)} 
\Hom_{\mb C} (\pi^*, \mathbb C) \otimes \mu^{\epsilon} \otimes {\nu}^{t}$$
must be contained in the Harish-Chandra module
$$V(\Ind_{\w{Sp}(p, \mb R)\w{GL}(q-p)N}^{\w{Sp}(q, \mb R)} \pi \otimes \mu^{\epsilon} \otimes {\nu}^{t}).$$
Our Lemma is proved. $\Box$
\begin{thm}~\label{id} The invariant tensor product $I^{\infty}(\epsilon, t)_{\w{U}(q)} \otimes_{\w{Sp}(p, \mb R)} V(\pi)$ can be identified with
$$V=V(\Ind_{\w{Sp}(p, \mb R)\w{GL}(q-p)N}^{\w{Sp}(q, \mb R)} \pi \otimes \mu^{\epsilon} \otimes {\nu}^{t}).$$
\end{thm}
Proof:  Recall that  $I^{\infty}(\epsilon, t)$ contains  the compactly supported $C^{\infty}_{c}(\mc M_{\epsilon, t}, d [g_1] d [k_2])$. In view of Equ. ~\ref{mt}, it is straight forward to show that
$$C^{\infty}_{c}(\mc M_{\epsilon, t}, d [g_1] d [k_2]) \ot V(\pi)$$ is dense and smooth in the Hilbert space
$\Ind_{\w{Sp}(p, \mb R)\w{GL}(q-p)N}^{\w{Sp}(q, \mb R)} \pi \otimes \mu^{\epsilon} \otimes {\nu}^{t}.$
Fix a $\tilde U(q)$-type $\sigma$. By orthogonality,
$$C^{\infty}_{c}(\mc M_{\epsilon, t}, d [g_1] d [k_2])_{\sigma} \ot V(\pi)$$ must be dense in
$$(\Ind_{\w{Sp}(p, \mb R)\w{GL}(q-p)N}^{\w{Sp}(q, \mb R)} \pi \otimes \mu^{\epsilon} \otimes {\nu}^{t})_{\sigma},$$
which is finite dimensional. Hence 
$$C^{\infty}_{c}(\mc M_{\epsilon, t}, d [g_1] d [k_2])_{\sigma} \ot V(\pi)= V_{\sigma}.$$
It follows that
$$I^{\infty}(\epsilon, t)_{\w{U}(q)} \otimes_{\w{Sp}(p, \mb R)} V(\pi) \supseteq C^{\infty}_{c}(\mc M_{\epsilon, t}, d [g_1] d [k_2])_{\w{U}(q)} \ot V(\pi)= V.$$
By the previous lemma, 
$I^{\infty}(\epsilon, t)_{\w{U}(q)} \otimes_{\w{Sp}(p, \mb R)} V(\pi) = V.$ $\Box$\\
\\
I shall now make a final remark. Fix an irreducible $\tilde U(q)$ representation $\sigma$. Suppose that $f \in I^{\infty}(\epsilon, t)_{\sigma}$. Then $f \ot v$ can be constructed as follows. Let $\phi$ be an arbitrary function in the $\sigma$-isotypic subspace
$$(\Ind_{\w{Sp}(p, \mb R)\w{GL}(q-p)N}^{\w{Sp}(q, \mb R)} \pi \otimes \mu^{\epsilon} \otimes {\nu}^{t})_{\sigma}.$$
Then the value of $\phi$ on $\tilde U(q)$ must be in a finite dimensional subspace in $V(\pi)$. Let $\{ u_i \}$ be an orthonormal basis of this finite dimensional subspace. Then one can define
$$(f \ot v)(k)=  \sum_{i} ((f \ot v)(k), u_i) u_i \qquad (k \in \tilde U(q)).$$

\section{Induced Intertwining Operator and Positivity  }
Let $t \in \mathbb R$.
Recall that there is an intertwining operator 
$$A(\epsilon,t): I^{\infty}(\epsilon, t) \rightarrow I^{\infty}(\epsilon, -t),$$
with poles on reducible $t$ (see ~\cite{vw}, ~\cite{sahi}, ~\cite{boo}).
 If $I(\epsilon,t)$ is irreducible, $A(\epsilon, t)$ one-to-one and onto. For complementary series $C(\epsilon, t)$, $A(\epsilon, t)$ is positive definite, one-to-one and onto. Let $(\, , \,)$ be the natural sesquilinear pairing between
$I(\epsilon, t)$ and $I(\epsilon, -t)$. The inner product of $C(\epsilon, t)$ is given by the completion of
$$(f_1, f_2)_{\epsilon, t}=(f_1, A(\epsilon, t) f_2).$$
Now $A(\epsilon,t)$ intertwines
the actions of $\w{Sp}(p, \mb R)$ and $\w{Sp}(q, \mb R)$. We can extend it trivially and obtain
$$A(\epsilon, t): I^{\infty}(\epsilon, t) \otimes V(\pi) \rightarrow I^{\infty}(\epsilon, -t) \otimes V(\pi).$$
Let $\pi \in \Pi_u(\w{Sp}(p, \mb R))$.
Suppose that $\pi$ and $\mu^{\epsilon}$ are compatible and $I^{\infty}(\epsilon, \pm t) \ot V(\pi)$ are well-defined.
We can now construct an induced intertwining operator
$$A(\epsilon, t, \pi): V(\Ind_{\w{Sp}(p, \mb R)\w{GL}(q-p)N}^{\w{Sp}(q, \mb R)} \pi \otimes \mu^{\epsilon} \otimes {\nu}^{t}) \rightarrow 
V(\Ind_{\w{Sp}(p, \mb R)\w{GL}(q-p)N}^{\w{Sp}(q, \mb R)} \pi \otimes \mu^{\epsilon} \otimes {\nu}^{-t})$$
via the following commutative diagram
\begin{equation}~\label{cm}
\begin{CD}
  I^{\infty}(\epsilon, t)_{\tilde U(q)} \otimes V(\pi)  @>{A(\epsilon, t) } >> I^{\infty}(\epsilon, -t)_{\tilde U(q)} \otimes V(\pi) \\
@VV{\ot}V @VV{ \ot}V \\
 V(\Ind_{\w{Sp}(p, \mb R)\w{GL}(q-p)N}^{\w{Sp}(q, \mb R)} \pi \otimes \mu^{\epsilon} \otimes {\nu}^{t}) @>{A(\epsilon, t, \pi)} >> 
 V(\Ind_{\w{Sp}(p, \mb R)\w{GL}(q-p)N}^{\w{Sp}(q, \mb R)} \pi \otimes \mu^{\epsilon} \otimes {\nu}^{-t})
\end{CD}
\end{equation}
We call $A(\epsilon, t, \pi)$ the {\it induced intertwining operator}. Presumably, $A(\epsilon, t, \pi)$ is contained in one of the intertwining operators constructed by Knapp-Stein and Vogan-Wallach (~\cite{ks3}, ~\cite{vw}).
\begin{thm}~\label{wd1} Let $\pi \in \Pi_u(\w{Sp}(p, \mb R))$.
Suppose that $\pi$ and $\mu^{\epsilon}$ are compatible. Suppose that  $I^{\infty}(\epsilon, \pm t) \ot V(\pi)$ is well-defined.
Then the induced intertwining operator $A(\epsilon, t, \pi)$ is well-defined. $A(\epsilon, t, \pi)$ intertwines the actions of $\f {sp}(q, \mb R)$ and $\tilde U(q)$.
\end{thm}
Proof: To show that $A(\epsilon, t, \sigma)$ is well-defined, it suffices to show that $\sum_{i} f_i \ot v_i=0$ implies $\sum_{i} A(\epsilon, t) f_i \ot v_i=0$. Fix a $\sigma \in \Pi_u(\tilde U(q))$. Let $f_1, f_2$ be in $I^{\infty}(\epsilon, t)_{\sigma}$, using the mixed model.  Let $v_1, v_2 \in V(\sigma)$. Let $\{ u_i\}$ be an orthonormal basis as in the final remark of Section 4. Then the natural sesqulinear pairing
\begin{equation}
\begin{split}
 & (f_1 \ot v_1, (A(\epsilon, t)f_2) \ot v_2) \\
 = & \int_{\tilde U(q)/\tilde U(p) \tilde O(q-p)} \sum_{\ u_i } \int_{\w{Sp}(p, \mb R)/\mc C} f_1( k h_1) (\pi(h_1) v_1, u_i) d [h_1] 
 \int_{\w{Sp}(p, \mb R)/\mc C} \overline{A(\epsilon, t)f_2(k h_2)(\pi(h_2) v_2, u_i)} \\ 
 &  d [h_2] d [k] \qquad [Absolutely \ convergent] \\
 = & \int_{\tilde U(q)/\tilde U(p) \tilde O(q-p)} \int_{\w{Sp}(p, \mb R)/\mc C} \int_{\w{Sp}(p, \mb R)/\mc C} f_1(k h_1) (\pi(h_1) v_1, \pi(h_2) v_2) \overline{A(\epsilon, t) f_2(k h_2)} d [h_1] d [h_2] d [k] \\
 = & \int_{\tilde U(q)/\tilde U(p) \tilde O(q-p)} \int_{\w{Sp}(p, \mb R)/\mc C} \int_{\w{Sp}(p, \mb R)/\mc C} f_1(k h_2 h) (\pi(h) v_1, v_2) \overline{A(\epsilon, t)f_2 (k h_2)} d [ h] d [h_2 ] d [k] \\
 & [h=h_2^{-1} h_1] \\
 =& \int_{\w{Sp}(p, \mb R)/\mc C} \int_{\tilde U(q)/\tilde U(p) \tilde O(q-p)}  \int_{\w{Sp}(p, \mb R)/\mc C} f_1(k h_2 h)  \overline{A(\epsilon, t)f_2 (k h_2)} d [h_2 ] d [k] (\pi(h) v_1, v_2) d [h] \\
  & [by \, Fubini's \, Theorem] \\
 = & \int_{\w{Sp}(p, \mb R)/\mc C} (I(\epsilon, t)( h)f_1, A(\epsilon, t) f_2) (\pi(h) v_1, v_2) d [h] \\
 = & \int_{\w{Sp}(p, \mb R)/\mc C} (I(\epsilon, t)( h) A(\epsilon, t) f_1, f_2)(\pi(h) v_1, v_2) d [h] \\
 = & ( (A(\epsilon, t) f_1) \ot v_1, f_2 \ot v_2))
 \end{split}
 \end{equation}
 Therefore $\sum_{i} f_i \ot v_i=0$ implies $\sum_{i} A(\epsilon, t) f_i \ot v_i=0$. The converse is also true since $A(\epsilon, t)$ is a surjection onto $I^{\infty}(\epsilon, -t)_{\tilde U(q)}$.
 $A(\epsilon, t, \pi)$ is well-defined.  Notice that the Diagram (\!~\ref{cm}) commutes and the actions of $\f {sp}(q, \mb R)$ and $\tilde U(q)$ are preserved in this diagram. Therefore $A(\epsilon, t, \pi)$ intertwines the actions of $\f {sp}(q, \mb R)$ and $\tilde U(q)$. $\Box$ \\
 \\
 Obviously, $A(\epsilon, t, \pi)$ is onto. It must also be one-to-one by Diagram (\!~\ref{cm}).
\begin{lem}~\label{critical} Let $\pi \in \Pi_u(\w{Sp}(p, \mb R))$.
Suppose that $C(\epsilon, t_0)$ is in the complementary series and $I^{\infty}(\epsilon, \pm t_0) \ot V(\pi)$ is well-defined. Then $I^{\infty}(\epsilon, \pm t) \ot V(\pi)$ is well-defined for
all real $-|t_0| \leq t \leq |t_0|$. 
\begin{enumerate}
\item the induced intertwining operator $A(\epsilon, t, \pi)$ is continuous in the interval
$t \in [-|t_0|, |t_0|]$ (for each $\tilde U(q)$-type);
\item $A(\epsilon, t, \pi)$ is one-to-one and onto for any $t \in [-|t_0|, |t_0|]$;
\item $A(\epsilon, 0, \pi)$ is the identity.
\end{enumerate}
\end{lem}
The assumption that $C(\epsilon, t)$ is in the complementary series is not absolutely necessary. Our assertions remain true except at those $t$ where $I(\epsilon, t)$ is reducible. So $A(\epsilon, t, \pi)$ may be discontinuous and not one-to-one at reducible $t$.

\begin{thm}~\label{main0} Fix $\epsilon \in [0, 1)$. Let $\pi \in \Pi_u(\w{Sp}(p, \mb R))$ be compatible with $\mu^{\epsilon}$.
Suppose that $|t_0| < \frac{1}{2}-|\frac{1}{2}-|2 \epsilon-1|| $ if $q-p$ odd and $|t_0| < |\frac{1}{2}-|2 \epsilon-1||$ if $q-p$ even.   Suppose $I^{\infty}(\epsilon, \pm t_0) \ot V(\pi)$ is well-defined.
Then $A(\epsilon, t, \pi)$ is positive definite for all $|t| \leq |t_0|$. Hence
$\Ind_{\w{Sp}(p, \mb R)\w{GL}(q-p)N}^{\w{Sp}(q, \mb R)} \pi \otimes \mu^{\epsilon} \otimes {\nu}^{t}$ is unitarizable for all $|t| \leq |t_0|$. 
\end{thm}
Proof: The positivity of $A(\epsilon, t, \pi)$  can be established by a deformation argument. Notice that $I(\epsilon, t)$ remains irreducible for all $|t| \leq |t_0|$. \\
\\
First $A(\epsilon, t, \pi)$ is well-defined for all $|t| \leq |t_0|$. It is continuous with respect to $t$ and it is always one-to-one and onto.
Second, since $A(\epsilon, t, \pi)$ intertwines the $\tilde U(q)$ actions, $A(\epsilon, t, \pi)$ restricted onto an $\tilde U(q)$-type is onto and one-to-one, therefore nondegenerate. Thirdly
the signature of $A(\epsilon, t, \pi)$ restricted to each $\tilde U(q)$-type (finite dimensional) must remain the same  for $t \in [-|t_0|, |t_0|]$.
Finally $A(\epsilon, 0, \pi)$ restricted to each $\tilde U(q)$-type is the identity. \\
\\
So $A(\epsilon, t, \pi)$ restricted onto each $\tilde U(q)$-type is positive definite. Hence it is positive definite on the Harish-Chandra module level. Therefore, the form $(* \, , *)_{t}$ defined as
$( * \, , A(\epsilon, t, \pi) *)$ is a positive definite invariant form on
$$V(\Ind_{\w{Sp}(p, \mb R)\w{GL}(q-p)N}^{\w{Sp}(q, \mb R)} \pi \otimes \mu^{\epsilon} \otimes {\nu}^{t})$$
Therefore $\Ind_{\w{Sp}(p, \mb R)\w{GL}(q-p)N}^{\w{Sp}(q, \mb R)} \pi \otimes \mu^{\epsilon} \otimes {\nu}^{t}$ is unitarizable.
$\Box$
\section{A General Construction of  $A(\epsilon, t, \pi)$}
It turns out that one can construct $A(\epsilon, t, \pi)$ in a much more general context. In particular, one can construct $A(\epsilon, t, \pi)$ for all unitary $\pi$ when $t$ is sufficiently negative. The idea is as follows. {\bf Fix  $\epsilon \in [0, 1)$ and $ t <0$}. Let $\pi \in \Pi_u(\w{Sp}(p, \mb R))$ be compatible. Consider the intertwining operator
$$A(\epsilon, t): C^{\infty}_{c}(\mc M_{\epsilon, t}, d [g_1] d [k_2])_{\tilde U(q)} \subseteq I^{\infty}(\epsilon, t)_{\tilde U(q)} \rightarrow I^{\infty}(\epsilon, -t)_{\tilde U(q)}.$$
The invariant tensor product $C^{\infty}_{c}(\mc M_{\epsilon, t}, d [g_1] d [k_2]) \ot V(\pi)$ is always well-defined (see Definition \!~\ref{itp}). In particular 
$$C^{\infty}_{c}(\mc M_{\epsilon, t}, d [g_1] d [k_2])_{\tilde U(q)} \ot V(\pi)=V(\Ind_{\w{Sp}(p, \mb R)\w{GL}(q-p)N}^{\w{Sp}(q, \mb R)} \pi \otimes \mu^{\epsilon} \otimes {\nu}^{t}).$$
Now suppose that $I^{\infty}(\epsilon, -t)_{\tilde U(q)} \ot V(\pi)$ is well-defined. Then we can define the {\it induced intertwining operator} $A(\epsilon, t, \pi)$ by the following diagram:
\begin{equation}~\label{cds}
\begin{CD}
 C^{\infty}_{c}(\mc M_{\epsilon, t}, d [g_1] d [k_2])_{\tilde U(q)} \otimes V(\pi)   @>{A(\epsilon, t) } >> I^{\infty}(\epsilon, -t)_{\tilde U(q)} \otimes V(\pi) \\
@VV{\ot}V @VV{ \ot}V \\
 V(\Ind_{\w{Sp}(p, \mb R)\w{GL}(q-p)N}^{\w{Sp}(q, \mb R)} \pi \otimes \mu^{\epsilon} \otimes {\nu}^{t}) @>{A(\epsilon, t, \pi)} >> 
 V( \Ind_{\w{Sp}(p, \mb R)\w{GL}(q-p)N}^{\w{Sp}(q, \mb R)} \pi \otimes \mu^{\epsilon} \otimes {\nu}^{-t})
\end{CD}
\end{equation}
Strictly speaking, we should have chosen a different notation for $A(\epsilon, t, \pi)$. As we shall show in Theorem ~\ref{criticals}, $A(\epsilon, t, \pi)$ we construct in this section will coincide with $A(\epsilon, t , \pi)$ we constructed in the last section if $I^{\infty}(\epsilon, t)_{\tilde U(q)} \ot V(\pi)$ is also well-defined. This allows us to view $A(\epsilon, t, \pi)$ here as a generalization of the construction from the last section.
\begin{thm}~\label{wd2} Let $\pi \in \Pi_u(\w{Sp}(p, \mb R))$ and $t \leq 0$.
Suppose that $\pi$ and $\mu^{\epsilon}$ are compatible. Suppose that  $I^{\infty}(\epsilon, - t) \ot V(\pi)$ is well-defined.
Then  $A(\epsilon, t, \pi)$ is well-defined. In addition, $A(\epsilon, t, \pi)$ intertwines the action of $\f {sp}(q, \mb R)$ and $\tilde U(q)$.
\end{thm}
Proof: Fix a $\sigma \in \Pi_u(\tilde U(q))$. Let $f_1, f_2$ be in $C^{\infty}_{c}(\mc M_{\epsilon, t}, d [g_1] d [k_2])_{\sigma}$, using the mixed model.  Let $v_1, v_2 \in V(\sigma)$. Let $\{ u_i\}$ be an orthonormal basis as in the final remark of Section 4. Then the natural pairing
$$(f_1 \ot v_1, (A(\epsilon, t)f_2) \ot v_2) =( (A(\epsilon, t) f_1) \ot v_1, f_2 \ot\ v_2).$$
Therefore $\sum_{i} f_i \ot v_i=0$ implies $\sum_{i} A(\epsilon, t) f_i \ot v_i=0$. The Diagram (~\ref{cds}) commutes.
$A(\epsilon, t, \pi)$ is well-defined. $\Box$\\
\\
Notice from Definition \!~\ref{itp}, if $t$ is sufficiently negative,  $I^{\infty}(\epsilon, - t) \ot V(\pi)$ is always well-defined. 
So $A(\epsilon, t, \pi)$ is well-defined for $t$ sufficiently negative with possible poles at reducible points. For $t$ close to zero,  $A(\epsilon, t, \pi)$ can perhaps be constructed by analytic continuation. It is an interesting problem to identify this intertwining operator in the standard construction
(~\cite{kn}). \\
\\
In the general context, $A(\epsilon, t, \pi)$ may have a kernel. However, if $I(\epsilon, t)$ is irreducible, $A(\epsilon, t, \pi)$ will be surjective. Therefore $A(\epsilon, t, \pi)$ must be one-to-one and nondegenerate. 

\begin{thm}~\label{criticals} Let $\pi \in \Pi_u(\w{Sp}(p, \mb R))$ and $t_0 \leq 0$.
Suppose that $I^{\infty}(\epsilon, - t_0) \ot V(\pi)$ is well-defined. Then  
\begin{enumerate}
\item $I^{\infty}(\epsilon, - t) \ot V(\pi)$ is well-defined for
all real $ t \leq t_0$.
\item the induced intertwining operator $A(\epsilon, t, \pi)$ is continuous on the interval
 $(-\infty, t_0)$ except at those $t$ for which $I(\epsilon, t)$ is reducible;
\item $A(\epsilon, t, \pi)$ is onto and one-to-one for  $t \in (-\infty, t_0)$ except at those $t$ for which $I(\epsilon, t)$ is reducible.
\end{enumerate}
\end{thm}
The reader shall compare this Lemma with Lemma ~\ref{critical}. \\
\\
Proof: The first assertion follows from Definition ~\ref{itp}. The second assertion follows from the computation in Theorems ~\ref{wd1}
and ~\ref{wd2}. We shall now prove the third assertion. Fix a $\sigma \in \Pi_u(\tilde U(q))$. Notice than $$
A(\epsilon, t) C^{\infty}_{c}(\mc M_{\epsilon, t}, d [g_1] d [k_2])_{\sigma} \subseteq I^{\infty}(\epsilon, -t)_{\sigma}.$$
Hence \begin{equation}
\begin{split}
 & A(\epsilon, t) C^{\infty}_{c}(\mc M_{\epsilon, t}, d [g_1] d [k_2])_{\sigma} \ot V(\pi) \\
 \subseteq  & I^{\infty}(\epsilon, -t)_{\sigma} \ot V(\pi) \\
 = & [\Ind_{\w{Sp}(p, \mb R)\w{GL}(q-p)N}^{\w{Sp}(q, \mb R)} \pi \otimes \mu^{\epsilon} \otimes {\nu}^{-t}]_{\sigma} \\
 = & C^{\infty}_{c}(\mc M_{\epsilon, -t}, d [g_1] d [k_2])_{\sigma} \ot V(\pi).
 \end{split}
 \end{equation}
  It suffices to prove that 
$$A(\epsilon, t) C^{\infty}_{c}(\mc M_{\epsilon, t}, d [g_1] d [k_2])_{\sigma} \ot V(\pi)=C^{\infty}_{c}(\mc M_{\epsilon, -t}, d [g_1] d [k_2])_{\sigma} \ot V(\pi).$$
We will show that every element on the right hand side can be approximated by elements on the left hand side. Since both sides are finite dimensional, they must be equal. \\
\\
Now let $f \in C^{\infty}_{c}(\mc M_{\epsilon, -t}, d [g_1] d [k_2])_{\sigma} \subset I^{\infty}(\epsilon, -t)_{\sigma}$. Then $A(\epsilon, t)^{-1} f \in I^{\infty}(\epsilon, t)_{\sigma}$. One can choose a sequence $\phi_i$  in $C^{\infty}_{c}(\mc M_{\epsilon, t}, d [g_1] d [k_2])_{\sigma}$ such that $\phi_i \rightarrow A(\epsilon, t)^{-1} f$ under the Frechet topology in $I^{\infty}(\epsilon, t)$. Then
$$A(\epsilon, t) \phi_i \rightarrow f$$
under the Frechet topology in $I^{\infty}(\epsilon, -t)$ (see ~\cite{vw}). In particular, $A(\epsilon, t) \phi_i$ converges to $f$ uniformly in the compact picture. In the mixed model 
$$ \| (A(\epsilon, t) \phi_i- f)(k_2 g_1) (\det g_1 g_1^t +1)^{\frac{n+1+t}{4}} \|_{\sup} \rightarrow 0.$$
Now it is easy to see that $A(\epsilon, t) \phi_i \ot v \rightarrow f \ot v$ uniformly on $\w{U}(q)$.
Hence $A(\epsilon, t) \phi_i \ot v \rightarrow f \ot v$ in the Hilbert space
$$\Ind_{\w{Sp}(p, \mb R)\w{GL}(q-p)N}^{\w{Sp}(q, \mb R)} \pi \otimes \mu^{\epsilon} \otimes {\nu}^{-t}.$$
$\Box$\\
\\
Observe that
$$A(\epsilon, t) C^{\infty}_{c}(\mc M_{\epsilon, t}, d [g_1] d [k_2]) \subseteq I^{\infty}(\epsilon, -t), \qquad
C^{\infty}_{c}(\mc M_{\epsilon, t}, d [g_1] d [k_2]) \subseteq I^{\infty}(\epsilon, t).$$
When $I^{\infty}(\epsilon, t)_{\tilde U(q)} \ot V(\pi)$ is also well-defined, the $A(\epsilon, t, \pi)$ here coincides with the $A(\epsilon, t, \pi)$ constructed in the last section. So there is no ambiguity. 
\begin{thm}~\label{main} Let $\epsilon \in [0, 1)$. Let $\pi \in \Pi_u(\w{Sp}(p, \mb R))$ be compatible with $\mu^{\epsilon}$.  Suppose that $0 \geq t > t_0$ where $t_0=-\frac{1}{2}+|\frac{1}{2}-|2 \epsilon-1|| $ if $q-p$ is odd and $t_0=-|\frac{1}{2}-|2 \epsilon-1||$ if $q-p$ is even. Suppose $I^{\infty}(\epsilon, 0) \ot V(\pi)$ is well-defined.
Then $A(\epsilon, t, \pi)$ is positive definite. Hence
$\Ind_{\w{Sp}(p, \mb R)\w{GL}(q-p)N}^{\w{Sp}(q, \mb R)} \pi \otimes \mu^{\epsilon} \otimes {\nu}^{t}$ is unitarizable. 
\end{thm}
Proof: Under our assumption, $I(\epsilon, t)$ is irreducible and unitarizable. By Theorem ~\ref{criticals}, $I^{\infty}(\epsilon, 0) \ot V(\pi)$ is well-defined implies that $I^{\infty}(\epsilon, -t) \ot V(\pi)$ is well-defined.
This theorem can be proved essentially the same way as Theorem \!~\ref{main0}. $\Box$ \\
\\
I shall make a final remark. The assumptions in this theorem are considerably weaker than the assumptions in Theorem ~\ref{main0}. Over the boundary point $t_0$ for which $I(\epsilon, t_0)$ is reducible, if one can choose a continuous $A(\epsilon, t)$ such that $A(\epsilon, t_0, \pi)$ is surjective, our unitarity theorem can be carried over to the next interval $(t_1, t_0)$. Here $I(\epsilon, t)$
must remain irreducible on $(t_1, t_0)$.
\section{Induced Complementary Series}
Now we can apply Theorem ~\ref{main} to build induced complementary series. Let us first give some results concerning
complementary principal series. 
\begin{thm}
Let $\sigma=\{\sigma_i\}_{i=1}^{n}$ be a compatible set of characters of the component group of $\w{GL}(1)$, namely 
$\exp 4 \pi \sqrt{-1} \sigma_i m=  \exp 4 \pi \sqrt{-1} \sigma_j  m  \ (m \in \mb Z \cong \mc C )$. Put $t_1= | 2 \sigma_1-1|$ and $t_i=\frac{1}{2}-|\frac{1}{2}-|2 \sigma_i-1||$ for $i >1$.
Then $$I(\sigma, v)=\Ind_{\w{GL}(1)^{n} N }^{\w{Sp}(n, \mb R)} \otimes ({\mu}^{\sigma_i} \otimes {\nu}^{v_i})$$ is unitarizable if
$v_i \in (-t_i, t_i)$ for every $i$.
\end{thm}
Proof: We prove this by induction on $n$. First of all, if $n=1$, $\Ind_{\w{GL}(1)}^{\w{Sp}(1, \mb R)} {\mu}^{\sigma_1} \otimes {\nu}^{v_1}$ is unitary if and only if $v_1 \in [-t_1, t_1]$, by a result due to Puk\'anszky (~\cite{puk}). Now suppose our theorem is true for $n-1$.
Let $\pi=\Ind_{\w{GL}(1)^{n-1} N }^{\w{Sp}(n-1, \mb R)} \otimes_i ({\mu}^{\sigma_i} \otimes {\nu}^{v_i})$ with $v_i \in (-t_i, t_i)$. Then
 leading exponents $\lambda$ of $\pi$ must be of the form $ w(v)- \rho(n-1)$, with $w$ a Weyl group element and $\rho(n-1)$ half sum of positive roots.
 Hence $\Re(\lambda) \prec -(n-2, n-3, \ldots, 1, 0)$. Now consider $I(\sigma_n, t; 2n-1)$ of $\w{Sp}(2n-1, \mb R)$ with $t \leq 0$. By Sahi's Theorem, 
 $I(\sigma_n, t; 2n-1)$ is unitarizable if $t \in (-t_n, t_n)$.
 Clearly,
 $$\Re(\lambda) -\mathbf{ \frac{2n-1+1}{2}}+ 2 \rho(n-1) \prec 0$$
 Hence $I(\epsilon_n, 0; 2n-1) \otimes_{\w{Sp}(n-1, \mb R)} \pi$ is well-defined. By Theorem ~\ref{main},
 $$\Ind_{\w{GL}(1)^{n} N }^{\w{Sp}(n, \mb R)} \otimes_{i=1}^n ({\mu}^{\sigma_i} \otimes  {\nu}^{v_i})=\Ind_{\w{GL}(1) \w{Sp}(n-1, \mb R) N}^{\w{Sp}(n, \mb R)} {\mu}^{\sigma_n} \otimes {\nu}^{v_n} \otimes \pi$$
 is unitarizable if $v_n \in (-t_n, t_n)$. $\Box$ \\
 \\
 For $\sigma_i=0, \frac{1}{2}$, the representation $I(\sigma, v)$ becomes a representation of the linear group. In this case, our theorem says nothing about the complementary series while there are plenty of complementary series at least when all $\sigma_i=0$ (~\cite{ko}).  For $\sigma_i=\frac{1}{4}, \frac{3}{4}$, our theorem says that $I(\sigma, v)$ is unitarizable if $v \in (-\frac{1}{2}, \frac{1}{2})^n$.
 These are the genuine principal complementary series of $Mp(n, \mb R)$. In particular, complementary series of size $(-\frac{1}{2}, \frac{1}{2})^n$ exist in every genuine principal series. Some of these complementary series were discussed in ~\cite{abptv} and ~\cite{barba}. \\
 \\
 Now let us deal with degenerate principal series. The theorem we are about to state is not the most general one. Our goal here is to show how one can build complementary series inductively.  
 \begin{thm}~\label{mu1}
  Let $L=\prod_{i=1}^{l} GL(r_i) Sp(r_0, \mb R)$ with 
 $$r_1 \leq r_2 \leq \ldots \leq r_l, \qquad \sum_{i=0}^{l} r_i=n.$$
 Let $\sigma_0$ be a tempered irreducible representation of  $\w{Sp}(r_0, \mb R)$. 
 Let the character ${\mu}^{\epsilon_i}$ of $\w{GL}(r_i)$ be compatible with $\sigma_0$.  Let 
 $$I(\epsilon, v, \sigma_0)=\Ind_{\w{Sp}(r_0, \mb R) \Pi_{i=1}^{l} \w{GL}(r_i) N}^{\w{Sp}(n, \mb R)} \sigma_0 \otimes (\otimes_{i=1}^{l} {\mu_i}^{\epsilon_i} \otimes  {\nu_i}^{v_i}) .$$
 Let $t_i=\frac{1}{2}-|\frac{1}{2}-|2 \epsilon_i-1||$ if $r_i$ is odd, $t_i=|\frac{1}{2}-|2 \epsilon_i-1||$ if $r_i$ is even. Then $I(\epsilon, v, \sigma_0)$ is unitarizable if $|v_i| < t_i$ all every $i \in [1, l]$.
 \end{thm}
 Proof: We proceed by induction on $i$. When $i=1$, by Lemma ~\ref{itp1}, $I(\epsilon_1, v_1; 2r_0+r_1) \otimes_{\w{Sp}(n_0, \mb R)} V(\sigma_0)$ is well-defined for $|v_1| \leq t_1$. By Theorem ~\ref{main0}, $I(\epsilon_1, v_1, \sigma_0)$ is unitarizable. In addition, leading exponents of $I(\epsilon_1, v_1, \sigma_0)$ are bounded by
 $$(\frac{r_1-1}{2}+|v_1|, \frac{r_1-1}{2}-|v_1|, \frac{r_1-3}{2}+|v_1|, \frac{r_1-3}{2}-|v_1| \ldots, 0, \ldots, 0) - \rho(r_0+r_1) $$
 $$ \preceq ( -r_0-\frac{r_1}{2}, -r_0-\frac{r_1}{2}+1, \ldots, 0)$$
 Here the last sequence has increment $1$. 
 It can be easily checked that
 $$( -r_0-\frac{r_1}{2}, -r_0-\frac{r_1-2}{2}, \ldots, 0)-\mathbf{\frac{2r_0+2r_1+r_2+1}{2}} \prec -2 \rho(r_0+r_1)$$
 since $r_2 \geq r_1$. So 
 $$I(\epsilon_2, 0; 2r_0+2 r_1+r_2) \otimes_{\w{Sp}(r_0+r_1, \mb R)} V(I(\epsilon_1, v_1, \sigma_0))$$
 is well-defined. By Theorem ~\ref{main}, $I(\epsilon_1, \epsilon_2, v_1, v_2, \sigma_0)$ is unitary. Our theorem follows by continuing this process. $\Box$
 
\commentout{

}


\begin{thebibliography}{99} 
\bibitem{abptv}[ABPTV] J. Adam, D. Barbasch, A. Paul, P. Trapa, D. Vogan \lq\lq Unitary Shimura correspondences for split real groups \!\rq\rq, {\it J. Amer. Math. Soc. } Vol. 20  2007, (701 - 751). 
\bibitem{barba}[B8] D. Barbasch \lq\lq Unitary Spherical Spectrum For Split
Classical Groups \!\rq\rq preprint, 2008.
\bibitem{bss}[BSS] D. Barbasch, S. Sahi, B. Speh \lq\lq Degenerate series representations for ${\rm GL}(2n,R)$ and Fourier analysis, \!\rq\rq  {\it Symposia Mathematica, Vol. XXXI }(Rome, 1988), Sympos. Math., XXXI, Academic Press, London, 1990, (45-69). 
\bibitem{bar}[Bar] V. Bargmann \lq\lq Irreducible unitary representations of the Lorentz group \rq\rq, {\it Annals of Math.}, (Vol 48), 1947 (568-640).
\bibitem{boo}[BOO] T. Branson, G. Olafsson, B. {\O}rsted, \lq\lq Spectrum Generating Operators and Intertwining Operators for Representations Induced from a Maximal Parabolic Subgroup, \!\rq\rq
{\it Journal of Functional Analysis}, (Vol. 135), 1996, (163-205).
\bibitem{du}[Du] M. Duflo, \lq\lq R\'epr\'esentation unitaires irr\'eductibles des groupes simples complexes de rang deux, \!\rq\rq   {\it Bull. Soc. Math. France}  107  (1979), no. 1, 55--96.

\bibitem{theta}[He00] Hongyu He, \lq\lq  Theta 
Correspondence I-Semistable Range: Construction and Irreducibility 
\rq\rq,  {\it Communications in Contemporary Mathematics } (Vol 2), 2000, 
(255-283).
\bibitem{unit}[Heu] Hongyu He, \lq\lq Unitary Representations and Theta Correspondence for
Classical Groups of Type I\rq\rq, {\it Journal of Functional Analysis}, Vol 199, Issue 1, (2003), 92-121.
\bibitem{nu}[Henu] Hongyu He {\it Unipotent Representations and Quantum Induction}, preprint, http://www.arxiv.org/math.GR/0210372, 2002.
\bibitem{comp}[Hec] Hongyu He, \lq\lq Restrictions of Degenerate Principal Series of the Universal Covering of the Symplectic Group \!\rq\rq, 2008, http://math.lsu.edu/~hongyu/pub/complementary.pdf (1-15).
\bibitem{hl}[HL] R. Howe, S. Lee, \lq\lq Degenerate principal series representations of ${\rm GL}_n(\mb C)$ and ${\rm GL}\sb n(\mb R)$,\rq\rq {\it  J. Funct. Anal. } 166  (1999),  no. 2, 244-309. 
\bibitem{howe}[Ho89] R. Howe, \lq \lq Transcending Classical Invariant 
Theory, \rq \rq
{\it Journal of Amer. Math. Soc.} (V. 2), 1989 (535-552).
\bibitem{john}[Joh] K. Johnson, \lq\lq Degenerate principal series and compact groups, \!\rq\rq {\it Math. Ann.}  287  (1990),  no. 4, 703-718.
\bibitem{kn}[KN] A. Knapp, {\it Representation Theory on Semisimple
Groups: An Overview Based on Examples} \, Princeton University Press, 
1986.
\bibitem{ks}[KS] A. Knapp, E. Stein, \lq\lq 
Intertwining operators for semisimple groups \!\rq\rq
{\it Ann. of Math.} (2) 93 (1971), 489--578. 
\bibitem{ks2}[KS2] A. Knapp, E. Stein, \lq\lq Intertwining operators for semisimple groups. II,\!\rq\rq {\it Invent. Math.}  60  (1980), no. 1, 9--84.
\bibitem{ks3}[KS3] A. Knapp, E. Stein, \lq\lq Some new intertwining operators for semisimple groups,\rq\rq  {\it Noncommutative harmonic analysis and Lie groups } (Marseille, 1980),  pp. 303-336, Lecture Notes in Math., 880, Springer, Berlin-New York, 1981. 
\bibitem{ko}[Kostant] B. Kostant \lq\lq On the existence and irreducibility of certain series of representations, \!\rq\rq {\it Bull. Amer. Math. Soc.} 75 (1969) 627-642.
\bibitem{kr}[KR] S. Kudla and S. Rallis, \lq\lq Degenerate Principal Series
and Invariant Distribution\rq\rq, {\it Israel Journal of Mathematics,} (Vol. 69, No. 1),
1990, (25-45).
\bibitem{lee}[Lee] S. T. Lee \lq\lq Degenerate principal series representations of ${\rm Sp}(2n,  R)$,\rq\rq {\it Compositio Math.}  103  (1996),  no. 2, (123-151).
\bibitem{lz}[LZ]  S. T. Lee, C-B. Zhu, \lq\lq Degenerate Principal Series and 
Local Theta
Correspondence II \rq\rq, {\it Israel Journal of Mathematics, 100}, 1997, (29-59).
\bibitem{li2}[LI89] J-S. Li, \lq \lq Singular Unitary Representation of 
Classical Groups, \rq \rq
{\it Inventiones Mathematicae} (V. 97), 1989  (237-255).
\bibitem{oz}[OZ] B. {\O}rsted, G. Zhang, \lq\lq Generalized Principal Series Representations and Tube Domain,\rq\rq {\it Duke Math. Journal}, (Vol. 78), 1995 (335-357).
\bibitem{puk}[Puk]  L. Puk\'anszky, \lq\lq The Plancherel formula for the universal covering group of ${\rm SL}(2, R),$ \rq\rq {\it Math. Ann. } 156  1964, (96-143).
\bibitem{sahi}[Sahi] S. Sahi, \lq\lq Unitary Representations on the Shilov Boundary of a Symmetric Tube Domain\rq\rq, {\it Representation Theory of Groups and Algebras}, (275-286), {\bf Comtemp. Math. 145}, Amer. Math. Soc., Providence, RI, 1993.
\bibitem{sv}[SV] B. Speh, D. Vogan, \lq\lq  Reducibility of generalized principal series representations, \rq\rq {\it  Acta Math.}  145  (1980), no. 3-4, 227-299.
\bibitem{vogan}[Vogan] D. Vogan, \lq\lq The Unitary Dual of $ GL(n)$ over an Archimedean
field \!\rq\rq {\it Invent. Math 83} Page 449-505 (1986).
\bibitem{vogang2}[Vog2] D. Vogan, \lq\lq The unitary dual of $G\sb 2$, \!\rq\rq {\it  Invent. Math.}  116  (1994),  no. 1-3, 677-791.
\bibitem{vw}[VW] D. Vogan and N. Wallach, \lq\lq Intertwining Operators for Real Reductive Groups, \!\rq\rq {\it Advances in Mathematics}, Vol 82, 1990 (203-243).
\end{thebibliography}
\end{document}